\documentclass[a4paper,10pt]{article} 
\usepackage{amsmath,amsthm,amsfonts,amssymb,color} 
\usepackage{graphicx} 
\usepackage[all]{xy}

 \relax
\input{amssym.def}

\newtheorem{defi}{Definition}[section]

\newtheorem{exem}[defi]{Example}
\newtheorem{exems}[defi]{Examples}
\newtheorem{lema}[defi]{Lemma}

\newtheorem{prop}[defi]{Proposition}
\newtheorem{teor}[defi]{Theorem}
\newtheorem{obs}[defi]{Remark}
\newtheorem{obss}[defi]{Remarks}

\newcommand{\huc}{\hspace*{.1cm}}

\newcommand{\n}{\neg} 

\newcommand{\su}{\mbox{$\huc\subseteq\huc$}}

\newcommand{\dm}{\noindent {\bf Proof:} }

\newcommand{\bc}{\begin{center}}
\newcommand{\ec}{\end{center}}

\newcommand{\IF}{{I\kern-0.3emF}}
\newcommand{\IK}{{I\kern-0.3emK}}

\newcommand{\IQ}{{I\kern-0.3emQ}}
\newcommand{\IR}{{I\kern-0.3emR}}

\newcommand{\lc}{{{\mathcal L}}}

\newcommand{\V}{\mathcal V}

\newcommand{\rar}{\rightarrow}

\newlength{\chico}

\begin{document}

\title{Matrix characterization of\\Ciuciura{'}s paraconsistent hierarchy 
$\textsf{Ciu}^n$}
\author{V\'\i ctor Fern\'andez; Gabriela Eisenberg\\
\small Basic Sciences Institute (Mathematical Area)\\ 
\small Av. Ignacio de la Roza (Oeste) 230, CP: 5400\\ \small San Juan, Argentina\\
\small E-mail:~{vlfernan@ffha.unsj.edu.ar}; {gabriela.eisenberg@gmail.com}\\}
\date{\today}

\maketitle

\begin{abstract}

\noindent In this paper, we will prove that the logics of the family $\textsf{Ciu}^n$:=$\{Ciu^n\}_{n \in \omega}$ of paraconsistent Ciuciura{'}s Logics (defined in  \cite{ciu:20} by means of bivaluations) can be alternatively defined by means of finite matrices. This result arises from the characterization of the truth-values of the involved matrices (relative to each $Ciu^n$-logic) as being specific finite sequences of elements of the set $2$ := $\{0,1\}$. Moreover, we will show along the paper that this characterization is related to the well-known standard Fibonacci Sequence, which is presented here by means of its binary expansion.  \end{abstract}

\normalsize 

\section{Introduction and Preliminaries}  \label{preliminares}

Many well-known paraconsistent logics in the literature are characterized by the fact that they do not validate any of the following principles:

\noindent ${\bf (NCP)}$ $\n (\varphi \wedge \n \varphi)$ (Non-contradiction Principle)

\noindent ${\bf (EP)}$ $\n \varphi \wedge \varphi \supset \psi$ (Principle of Explosion/Trivialization)

This fact determines that they also fail to validate the Principle of {\it Expansion of Double Negation}:

\noindent ${\bf (DNE)}$ $\varphi \supset \n \n \varphi$

Among these logics, we mention da Costa's hierarchy $\{C_n\}_{1 \leq n \leq \omega}$ (see \cite{dac:74}), Sette's logic $P^1$, given in 
\cite{set:73}, and also a generalization of it, denoted by $P^n$ (or by  $I^0 P^n$; see 
\cite{fer:22}, 
\cite{fer:con:03}). 

With a different perspective, J. Ciuciura defined in \cite{ciu:20} a family of logics (denoted in this paper by $\textsf{Ciu}^n$, with $n \geq 1$), motivated by the following approach: 
the essential principle that should define paraconsistency is the non-validation of
${\bf (EP)}$. Actually, the common language related with the $Ciu^n$-logics is simply $L(C)$, based on $C$ = $\{\n, \supset\}$, which does not possess a connective expressing the conjunction, necessary for the statement of ${\bf (NCP)}$. On the other hand (and denoting the consequence relation of $Ciu^n$ by $\models_{Ciu^n}$), we have that ${\bf (EP)}$ can be expressed as $\n \varphi, \varphi \models_{Ciu^n} \psi$, without necessity of  $\wedge$. In a similar way, 
${\bf (DNE)}$ can be expressed by $\varphi \models_{Ciu^n} \n \n \varphi$. And none of these consequences are valid in $Ciu^n$ with $n \geq 1$.

We mention here the following fact: actually, the consequence relation $\models_{Ciu^n}$ was defined in \cite{ciu:20} by means of two methods that are equivalent ones (the proof of this fact obbeys to the usual Soundess/Completeness Theorems). These definions are based respectively on:

$(a)$ A Hilbert-Axiomatics.

$(b)$ A {\it semantics of bivaluations}. 

With respect to the latter notion, it should be understood that it is refers to a family $S_n$, of maps from the set $L(C)$ (of formulas) to the set $2$ = $\{0,1\}$. These functions, {\it that are not necessarily homomorphics}, determine every relation $\models_{Ciu^n}$. 

Let us note here that, even when the $Ciu^n$-logics are very well motivated, and even when their definitions are natural enough, it still remains open the following technical problem: {\it Are these logics characterized by means of finite matrices?} This problem is interesting in the context of the study of  paraconsistent logics because there is not  a general result that can be applied to all of them. For instance, none logic of the paraconsistent da Costa's hierarchy 
$\{C_n \}_{1 \leq n \leq \omega}$ is able to be expressed by using finite matrices (see \cite{arr:75}). On the other hand, all the logics $I^0P^n$ are directly defined using matrix techniques.

In this paper, we will demonstrate the following result, specifically referred to the $Ciu^n$-logics: {\it all the logics} of the family 
$\textsf{Ciu}^{n}$:=$\{Ciu^n\}_{n \geq 0}$ are expressed by means of finite matrices\footnote{Is should be anticipated here that, extending the definition given in \cite{ciu:20}, in this paper 
the family $\textsf{Ciu}^n$ includes the logic $Ciu^0$, that will be identified with the classical logic $CL$. 
We will explain this, in a more rigourous way, later.}. Moreover, the cardinality of such matrices is directly related with the well-known  {\it  Fibonacci{'}s Sequence $Fb(n)$}, $1 \leq n$, as we shall see.
To obtain all these results, let us previously fix a minimum of notation and basic definitions.  With respect to the formal language of the involved logics and to matrix logics, and taking \cite{fon:16} as basis, let us recall the following concepts:

\begin{defi}\label{lenguaje}\rm{ We consider as starting point the set $\omega$:=$\{0,1, 2,\dots\}$, and fix a countable set $\V$ = $\{p_i\}_{i \in \omega}$, of {\it atomic formulas}. The {\it signature} to be used (that is, the set of primitive connectives), common to all the logics of the family $\textsf{Ciu}^{n}$, is $C$ := $\{\n, \supset\}$; the {\it language} of such logics will be defined as the absolutely free algebra, generated by $C$ over $\V$, and it will be denoted by  
$L(C)$. }
\end{defi}

We will use lowercase greek letters $\varphi, \psi, \theta...$ as metavariables ranging on $L(C)$; in particular, the  letter $\alpha$ will be reserved to indicate elements of $\V$ (all these symbolos well be sub/over-indexed, if necessary). Bearing this in mind, $\n^k \varphi$ abbreviates $\underbrace{\n (\n \dots \n (\varphi)\dots)}_{\textrm{$k$ times}}$; $\n^0 \varphi$ is simply $\varphi$.  The {\it sets} of formulas will be denoted by uppercase greek letters $\Gamma$, $\Delta$, etc.

\begin{defi}\label{matrices-gral}\rm{A $C$-matrix is a pair $M$ = $({\bf A}, D)$, where ${\bf A}$ is an algebra which is similar to $L(C)$, and $D \su A$, being $A$ the support of ${\bf A}$.}
\end{defi}

The operations of every $C$-matrix $M$ will be denoted using $\supset$ and $\n$ too, to avoid excessive notation. Despite this, if it were necessary (in some proof where different $C$-matrices are compared, for instance), the involved operations will be  indexed. Besides that, let us recall that every $C$-matrix defines a consequence relation on $L(C)$:

%
%
%
%
%

\begin{defi}\label{cons-matricial}\rm{Given a  $C$-matrix $M$ = $({\bf A},D)$, the consequence relation
$\models_M \su \wp(L(C)) \times L(C)$ is defined as follows: $\Gamma \models_M \varphi$ if and only if, for every 
\\
{\it $M$-valuation}	 (i.e., every homomorphism $w: L(C) \to A$) such that
$w(\Gamma) \su D$, it holds that $w(\varphi) \in D$.}
\end{defi}

We will use later the following definitions and results referred to the consequence relations defined by matrices (see \cite{fon:16}):

\begin{defi}\rm{Two $C$-matrices $M_i$ = $({\bf A_i}, D_i)$ are {\it isomorphic} if and only if there is a bijection $f:A_1 \to A_2$ such that:

\noindent $(a)$ $f$ is an algebraic isomorphism between ${\bf A_1}$ y ${\bf A_2}$.

\noindent $(b)$ For every $x \in A_1$, $x \in D_1$ if and only if $f(x) \in D_2$.
}
\end{defi}

\begin{prop}\label{iso-mat-misma-cons}\rm{If $M_1$ and $M_2$ are isomorphic, then for every $\Gamma \cup \{\varphi\} \su L(C)$ it holds that 
$\Gamma \models_{M_1}\varphi$ if and only if $f(\Gamma) \models_{M_2} f(\varphi)$.
}
\end{prop}

Let us recall additionally that, since $L(C)$ is an  absolutely free algebra, every $M$-valuation $w$ is univocally determined by the values $w(\alpha)$, considering all the formulas $\alpha \in \V$. At this point we remark that, for every $C$-matrix $M$, the pair  $(L(C),\models_M)$ is an {\it abstract (propositional) logic}, cf. \cite{bro:sus:73}, \cite{fon:16}. That is, $\models_M$ verifies the following properties: 

\

\noindent $(Ext)$ If $\varphi \in \Gamma$, then $\Gamma \models_M \varphi$. \hfill(Extensivity)

\noindent$(Mon)$ If $\Gamma \models_M \varphi$ \ and \ $\Gamma \su \Delta$, then $\Delta \models_M \varphi$. \hfill(Monotonicity)

\noindent $(Tran)$ If $\Gamma \models_M \varphi$ \ and \ $\Delta \models_M \gamma$ for every $\gamma \in \Gamma$, then $\Delta \models_M \varphi$. 

\hfill(Transitivity)

Aditionally, every matrix logic $(L(C), \models_M)$ verifies:

\noindent $(Str)$ If $\Gamma \models_M \varphi$ then, for every {\it substitution} (i.e. endomorphism)

\noindent $\sigma: L(C) \to L(C)$, it holds that $\sigma(\Gamma) \models_M \sigma(\varphi)$. \hfill(Structurality)

\noindent $(Fin)$ If $M$ has finite support $A$, then: if $\Gamma \models_M \varphi$  there is $\Gamma_0 \su \Gamma$, $\Gamma_0$ finite, such that
$\Gamma_0 \models_M \varphi$. \hfill(Finitariness)

\

Obviously, a consequence relation verifying $(Ext)$, $(Mon)$ and $(Tran)$ can be defined by other ways, and not necessarily by means of matrices. In this paper we are interested in a very simple  (and natural) way of defining such consequence relations, which consists of considering families of 
{\it non-homomorphic bivaluations}. This method can be formalized in the following way, adapting the approach of \cite{car:con:16}:

\begin{defi}\label{seman-bival}\rm{A {\it bivaluated semantics} is a non-empty set of functions

\noindent $S$ = $\{v_i: L(C) \to \{0,1\}\}_{i \in I}$. The elements of $S$ are called {\it bivaluations}, and they not need to be homomorphisms (in fact, we have not defined any algebraic structure in the set $2$ = $\{0, 1\}$, yet). The set $S$ defines the consequence relation $\models_S \su \wp(L(C)) \times L(C)$ as follows: for every $\Gamma \cup \{\varphi\}\su L(C)$, $\Gamma \models_S \varphi$ if and only if, for every $v \in S$ such that $v(\Gamma)$ = $\{1\}$, it holds that $v(\varphi)$ = $1$.}
\end{defi}

It is easy to prove that, for every bivaluated semantics $S$, $\models_S$ verifies $(Ext)$, $(Mon)$ and $(Tran)$ (and therefore 
$(L(C),\models_S)$ is an abstract logic, too). However, it is not warranted that $\models_S$ verifies $(Str)$ and/or $(Fin)$. The bivaluated semantics are very frequently used in the study of certain non-classical logics, and specially in Paraconsistent Logic. In fact, the well-known {\it semantics of quasi-matrices} (see \cite{dac:alv:77}) can be understood as a particular case of bivaluated semantics, for instance.

\

We conclude this brief introduction remarking that we will use often the boolean algebra ${\bf 2}$=$(2,\{+,\cdot,-,0,1\})$ (the meaning of the operations are the usual ones). In this algebra, it can be defined in the usual way  $x\rar y$:=$-x + y$. So, we avoid to use the same symbols for the connectives of  $C$ and for the operations of  ${\bf 2}$. This differentiation will be very convenient for the rest of the paper. Moreover, we will be focused on the  algebra 
${\bf 2}^{\{-,\rar\}}$,  the {\it $\{-,\rar\}$-reduct of ${\bf 2}$}. We will work specifically with it
because the $Ciu^n$-bivaluations 
can be, {\it eventually}, homomorphisms between $L(C)$ and ${\bf 2}^{\{-,\rar\}}$ (this makes sense now, considering that both algebras are similar). 
Using all the definitions and results given above, we will start the analysis of the hierarchy $\textsf{Ciu}^n$.

\section{On Ciuciura{'}s logics $Ciu^n$}

The $Ciu^n$-logics, with $n \in \omega$, can be defined by means of bivaluations in the following way
 (see \cite{ciu:20}):

\begin{defi}\label{ciun-bivaluaciones}\rm{Let $n \in \omega$, fixed. A {\it $Ciu^n$-bivaluation} is any map 
\\
$v: L(C) \to \{0,1\}$ that fulfill the following conditions, for every $\varphi$, $\psi \in L(C)$:

\noindent ${\bf (1)}$ If $v(\n \varphi)$ = $0$, then $v(\varphi)$ = $1$.

\noindent ${\bf (2.n)}$ If $v(\n^{n+1} \varphi)$ = $1$, then $v(\n^n\varphi)$ = $0$.

\noindent ${\bf (3)}$ If $v(\n(\varphi \supset \psi))$ = $1$, then $v(\varphi \supset \psi)$ = $0$.

\noindent ${\bf (4)}$ $v(\varphi \supset \psi)$ = $1$ if and only if $v(\varphi)$ = $0$ or $v(\psi)$ = $1$. 

\

The set of all the $Ciu^n$-bivaluations will be denoted by $S_n$. It determines the relation $\models_{S_n}$, according Definition 
\ref{seman-bival}.
}
\end{defi}

\begin{defi}\label{defi-logicas-ciu-n}\rm{For every $n \geq 0$, {\it Ciuciura's logic $Ciu^n$} is defined as being the pair $Ciu^n$:=$(L(C),\models_{S_n})$. The family of these logics will be indicated by $\textsf{Ciu}^n$. For the sake of brevity, every logic of this family will be mentioned usually as a ``$Ciu^n$-logic''.}
\end{defi}

\begin{obss}\label{ciu-0}\rm{With respect to the definition above:

\noindent $(a)$ In \cite{ciu:20}, J. Ciuciura considers just $n \geq 1$. We have included here $n$ = $0$. This inclusion does not generate conflicts in the results to be developed all along this paper. Moreover, it allows to identify the classical logic $CL$ with $Ciu^0$.

\noindent $(b)$ It is easy to see that for every $n$, $m \in \omega$, $n \leq m$ implies $S_n \su S_m$.

\noindent $(c)$ In the definition above, conditions ${\bf (1)}$ and ${\bf (2.n)}$ together imply that, for every  $k \geq n$, every $\alpha \in \V$, it holds that
$v(\n^k \alpha)$ = $0$ if and only if $v(\n^{k+1} \alpha)$ = $1$. This property is also valid for every $v \in S_n$, for every implication $\varphi \supset \psi$, because conditions  ${\bf (1)}$ and ${\bf (3)}$. 
Moreover, it holds for all the formulas of the form $\n^k \alpha$, with $k \geq n$, $\alpha \in \V$.
In an informal way we can say that all these formulas ``behave classically'' (with respect to the negation connective $\n$).

\noindent $(d)$ On the other hand, if $0 \leq k < n$ and $\alpha \in \V$, we can only warrant the following, from {\bf (1)}: if $v(\n^k \alpha)$ = $0$, then $v(\n^{k+1} \alpha)$ = $1$. So, if $k < n$, there may be $Ciu^n$-bivaluations verifying $v(\n^{k}\alpha)$ = $v(\n^{k+1}\alpha)$ = $1$.
%
%

\noindent $(e)$ From all the previously mentioned, note that every $Ciu^n$-bivaluation $v$ will be univocally determined by the values $v(\n^k \alpha)$, with 
$0 \leq k \leq n$, $\alpha \in \V$. Hovever, all these values {\it are not completely independent} from each other. 
This fact is, indeed, one of the main characteristics that will determine the connection  of the logics of the family $\textsf{Ciu}^n$ with Fibonacci{'}s sequence, as we shall see later. 

\noindent $(f)$ According to $(d)$, it is easy to see that the $Ciu^0$-bivaluations simply are the standard valuations (i.e. homomorphisms) of the classical logic 
$CL$. Then, $\models_{S_0}$ = $\models_{CL}$, as we were indicated in $(a)$.
}
\end{obss}

Let us see some examples of bivaluations, in different logics of the family $\textsf{Ciu}^{n}$.

\begin{exems} \label{bival-sencillas}
\noindent \rm{Let ${\bf 2}^{\{-,\rar\}}$ be the reduct indicated at the end of Section \ref{preliminares},
and consider $n$ = $1$. We define the map $v_1: L(C) \to 2$ by:

\

$v_1(\varphi)$ = 
$\left\{
\begin{array}{ll}
1 & \textrm{ if $\varphi \in \V$ }\\
1 & \textrm { if $\varphi$ = $\n \alpha$, with $\alpha \in \V$}\\
- \, v_1(\psi) & \textrm{ if $\varphi$ = $\n \psi$, with $\psi \notin \V$}\\
v_1(\psi) \rar  v_1(\theta) & \textrm{ if $\varphi$ = $\psi \supset \theta$}\\
\end{array}
\right.$

\

\

\noindent It is easy to check that $v_1 \in S_1$, and that $v_1$ is not an homomorphism from $L(C)$ to ${\bf 2}^{\{-,\rar\}}$. 

\noindent On the other hand, let us consider the map $v_2: L(C) \to 2$ given by:

\

$v_2(\varphi)$ = $\left\{
\begin{array}{ll}
1 & \textrm{ if $\varphi \in \V$ }\\
1 & \textrm { if $\varphi$ = $\n \alpha$, with $\alpha \in \V$}\\
{ 0} & \textrm {if $\varphi$ = $\n^2 \alpha$, with $\alpha \in \V$}\\
- v_2(\psi) & \textrm{ if $\varphi$ = $\n \psi$, with $\psi \neq \n^k \alpha$,} \\
& \textrm{$\alpha \in \V$, $k \in \{ 0,1\}$}\\
v_2(\psi) \rar v_2(\theta) & \textrm{ if $\varphi$ = $\psi \supset \theta$}\\
\end{array}
\right.$

\

\noindent Finally, let us define:

\

$v_3(\varphi)$ = $\left\{
\begin{array}{ll}
1 & \textrm{ if $\varphi \in \V$ }\\
1 & \textrm { if $\varphi$ = $\n \alpha$, with $\alpha \in \V$}\\
{ 1} & \textrm { if $\varphi$ = $\n^2 \psi$, with $\alpha \in \V$}\\
- v_3(\psi) & \textrm{ if $\varphi$ = $\n \psi$, with $\psi \neq \n^k \theta$,} \\
& \textrm{$\alpha \in \V$, $k\in \{0,1\}$}\\
v_3(\psi) \rar v_3(\theta) & \textrm{ if $\varphi$ = $\psi \supset \theta$}\\
\end{array}
\right.$ 

\

\noindent It can be checked that $\{v_1, v_2, v_3\} \su S_2$. On the other hand, $v_2 \notin S_1$, $v_3 \notin S_1$.
}
\end{exems}
It was proved in \cite{ciu:20} that $Ciu^1$ can be matricially characterized. Indeed, $Ciu^1$ is Sette{'}s Logic $P^1$ given in \cite{set:73}. In addition, $Ciu^0$ is also characterized by a matrix: since $Ciu^0$ = $CL$,
it is definable by means of the matrix $M_{CL}$ = $({\bf 2},\{1\})$ (the standard two-valued matrix of $CL$). In the next section we will give a matrix representation for $Ciu^2$, which will motivate the analyisis of all the $Ciu^n$-logics.

\section{A motivating example: matrix semantics for $Ciu^2$}\label{semantica-para-ciu2}

We will reformulate Definition \ref{seman-bival} now, applying it to the  specific case of the bivaluated semantics for $Ciu^2$:

\begin{defi}\label{semantica-s2}\rm{A {\it $Ciu^2$-bivaluation} is any map $v: L(C) \to \{0,1\}$ verifying the following conditions, for every  $\varphi$, $\psi \in L(C)$:

\noindent ${\bf (1)}$ If $v(\n \varphi)$ = $0$, then $v(\varphi)$ = $1$.

\noindent ${\bf (2{.}2)}$ If $v(\n^{3} \varphi)$ = $1$, then $v(\n^2\varphi)$ = $0$.

\noindent ${\bf (3)}$ If $v(\n(\varphi \supset \psi))$ = $1$, then $v(\varphi \supset \psi)$ = $0$.

\noindent ${\bf (4)}$ $v(\varphi \supset \psi)$ = $1$ if and only if $v(\varphi)$ = $0$ or $v(\psi)$ = $1$. 

\noindent The set of all the $Ciu^2$-bivaluations will be denoted by $S_2$, defining the consequence relation $\models_{S_2}$ as it was previously said.
}
\end{defi}

As it was previously remarked, ${\bf (2{.}2)}$ together with ${\bf (1)}$ imply that $v(\n^2 \varphi)$ = $0$ if and only if $v(\n^3 \varphi)$ = $1$.

We will give some examples of $Ciu^2$-bivaluations, in the sequel. For that, it is very convenient to establish the following definition:

\begin{defi}\label{k-n}{\rm 
For every $n \in \omega$, we define the  set of formulas $K_n^{\ast}$ by:

\noindent $K_n^{\ast}$:= $\{\varphi \in L(C): \varphi = \n^k \alpha$, $\alpha \in \V$, $0 \leq k < n\}$ \footnote{Note here that 
$K_0^{\ast}$ = $\emptyset$.}.
In the concrete case when $n = 2$,
$K_2^{\ast}$:=$\{\varphi \in L(C): \varphi = \n^k \alpha$, $\alpha \in \V$, $0 \leq k \leq 1\}$.
} 
\end{defi}

The sets $K_n^{\ast}$ defined above will allow us to vinculate everv ${Ciu}^n$-logic with Fibonacci{'}s  sequence, as we shall see. On the other hand, these sets can give us some simple examples of  ${\textsf{Ciu}}^n$-bivaluations. Let us see some of them (with $n$ = $2$) in the sequel. 

\noindent \begin{exem}\label{ejemplos-ciu2}\rm{
Let us see a couple of standard examples of non-homomorphical $Ciu^2$-bivaluations, following the same patterns of Example \ref{bival-sencillas}:

\

\noindent $(a)$ We define  $v_1 \in S_2$ as follows, considering $\varphi$ varying on $L(C)$: 

\noindent $v_1(\varphi) = \left\{
\begin{array}{ll}
1 & \textrm{ if $\varphi \in K_2^{\ast}$}\\
- v_1(\psi) & \textrm{ if $\varphi$ = $\n \psi$, $\psi \notin K_2^{\ast}$}\\
v_1(\psi)\rar v_1(\theta) & \textrm{ if $\varphi$ = $\psi \supset \theta$}\\
\end{array}
\right.$

\

\noindent $(b)$ In addition, we define $v_2 \in S_2$ in the following way:

\noindent $v_2(\varphi)$ = 
$\left\{
\begin{array}{ll}
0 & \textrm{ if $\varphi \in \V$ }\\
1 & \textrm { if $\varphi$ = $\n \alpha$, with $\alpha \in \V$}\\
1& \textrm{ if $\varphi$ = $\n \n \alpha$, with $\alpha \in \V$}\\
- v_2(\psi) & \textrm{ if $\varphi$ = $\n \psi$, with $\psi \notin \V$}\\
v_2(\psi) \rar v_2(\theta) & \textrm{ si $\varphi$ = $\psi \supset \theta$}\\
\end{array}
\right.$

\

\noindent Both functions are clearly well defined. Besides that, it can informally noted that they are $Ciu^2$-bivaluations reasoning as follows: according Definition \ref{semantica-s2}, an arbitrary map $v:L(C)\to  2$ must verify the following, in order to belong to $S_2$: first, for every formula $\varphi \notin K_2^{\ast}$ (``outside $K_2^{\ast}$''), $v$ should  ``behave as an homomorphism'' from $L(C)$ to ${\bf 2}^{\{-,\rar\}}$. 
In addition, for the elements of $K_2^{\ast}$, condition ${\bf (2{.}2)} $ only demands the following: if $\varphi \in K_2^{\ast}$, $v \in S_2$, {\it and additionally $v(\varphi)$ = $0$}, then $v(\n \varphi)$ = $1$. On the other hand, there is not any additional condition when $\varphi \in K_2^{\ast}$, $v(\varphi)$ = $1$. And all these requirementes are verified by $v_1$ and by  $v_2$.}
\end{exem}

We will develop the matrix semantics for $Ciu^2$, in the sequel. 

\begin{defi}\label{matriz-a2-defi}\rm{Consider the set $A_2$:=$\{\vec{x}=(x_0,x_1,x_2) \in 2^3: \vec{x} \textrm{ verifies $(\star)_2$}\}$ (here, condition $(\star)_2$ is: if $x_i$ = $0$, then $x_{i+1}$ = $1$, for $i$ = $0$, $1$). 
Additionally, the matrix $M_2$ is defined by $M_2$:=$({\bf A_2},D_2)$, with
\\
$D_2$:= $\{\vec{x}= (x_0,x_1,x_2)\in A_2: x_0 = 1\}$, and being ${\bf A_2}$:= $(A_2,\{\n,\supset\})$ the 
\\
$C$-algebra whose operations $\n$ and $\supset$ are defined as follows (for 
$\vec{x}$ = $(x_0,x_1,x_2)$, 
$\vec{y}$ = $(y_0,y_1,y_2) \in A_2$):

\noindent $(a)$ $\n \vec{x}$ =  $\n (x_0,x_1,x_2)$:=$(x_1,x_2,-x_2)$.

\noindent $(b)$ $\vec{x} \supset \vec{y}$ = $(x_0,x_1,x_2) \supset (y_0,y_1,y_2)$:=  $(x_0 \rar y_0, - (x_0 \rar y_0), x_0 \rar y_0)$
}
\end{defi}

\begin{obs}\label{M2-en-limpio}
\rm{Concerning the previous definition: 

\noindent $(a)$ The sets $A_2$ and $D_2$ are, respectively:

$\bullet$  $A_2$ = $\{(0,1,0); (0,1,1); (1,0,1); (1,1,0); (1,1,1)\}$.

$\bullet$ $D_2$ = $\{(1,0,1); (1,1,0); (1,1,1)\}$.

\noindent $(b)$ In addition, the operations of ${\bf A_2}$ are given in these truth-tables:

\bc
$\begin{array}{|c|c c c c c|} \hline
\n & (0,1,0) & (0,1,1) & (1,0,1) & (1,1,0) & (1,1,1)\\ \hline
& (1,0,1) & (1,1,0) & (0,1,0) & (1,0,1) & (1,1,0)\\
\hline
\end{array}
$

\

\

$\begin{array}{|c|c c c c c|} \hline
\supset & (0,1,0) & (0,1,1) & (1,0,1) & (1,1,0) & (1,1,1)\\ \hline
(0,1,0) & (1,0,1) & (1,0,1) & (1,0,1) & (1,0,1) & (1,0,1)\\
(0,1,1) & (1,0,1) & (1,0,1) & (1,0,1) & (1,0,1) & (1,0,1)\\
(1,0,1) & (0,1,0) & (0,1,0) & (1,0,1) & (1,0,1) & (1,0,1)\\
(1,1,0) & (0,1,0) & (0,1,0) & (1,0,1) & (1,0,1) & (1,0,1)\\
(1,1,1) & (0,1,0) & (0,1,0) & (1,0,1) & (1,0,1) & (1,0,1)\\
\hline
\end{array}
$
\ec

Let us note the following facts, taking into account Definition \ref{k-n}: 

\noindent $(c)$ For every $\varphi \notin K_2^{\ast}$, for every $M_2$-valuation $w$, it holds that $w(\varphi)$ = $(0,1,0)$ or $w(\varphi)$ = $(1,0,1)$. 

\noindent $(d)$ For every $\vec{x}$, $\vec{y} \in A_2$, $\vec{x} \supset \vec{y} \in D_2$ if and only if $\vec{x} \notin D_2$ or $\vec{y} \in D_2$.
}
\end{obs}

%
%
%
%
%
%
%

%
%
%
%
%
%
%
%
%
%
%

Let us prove that $\models_{S_2}$ = $\models_{M_2}$, starting from the following technical results:

\begin{lema}\label{lema-clave-ciu2-1}
\rm{Every bivaluation $v \in S_2$ determines a $M_2$-valuation $w_v$ such that, for every $\varphi \in L(C)$, $w_v(\varphi) \in D_2$
if and only if $v(\varphi)$ = $1$.}
\end{lema}

\dm Let $v$ be in $S_2$. For every $\alpha \in \V$, we define $w_v(\alpha)$:= $(v(\alpha), v(\n \alpha), v(\n^2 \alpha))$. It can easily seen that 
$w_v(\alpha) \in A_2$, because the conditions asked for $S_2$ in Definition \ref{semantica-s2}. 
That is, $w_v(\alpha)$ is well defined. In addition, let us extend homomorphically $w_v$ following the definitions of  $\n$ and $\supset$ in ${\bf A_2}$, obtaining so $w_v(\varphi)$ for every $\varphi \in L(C)$. We will prove now that, for every $\varphi \in L(C)$:
\bc $(\circledast)$ $w_v(\varphi) \in D_2$  if and only if $v(\varphi)$ = $1$. \ec
\noindent Our proof will be done by induction on the complexity of $\varphi$. From now on and all along the paper in this kind of proofs, IH abbreviates ``induction hypothesis''.
\\
\noindent \underline{Case 1}: $\varphi\in\mathcal{V}$. It holds by definition.
\\
\noindent \underline{Case 2}: $\varphi=\n\psi$ 

\noindent 2.1) Let us suppose first that $w_v(\varphi) \in D_2$. We have these possibilities:

\noindent 2.1.1) $w_v(\varphi)$=$(1,0,1)$. Thus, $w_v(\psi)$=$(1,1,0)$ or $w_v(\psi)$=$(0,1,0)$, cf. Remark \ref{M2-en-limpio}. Now:
\\
\noindent 2.1.1.1) Let us suppose $\psi\in\mathcal{V}$; if $w_v(\psi)$=$(1,1,0)$ then $v(\psi)=1$, by IH. From this, it can happen that
$v(\n\psi)=0$ or $v(\n\psi)=1$.
If $v(\varphi) =  v(\n\psi)=0$, then $v(\n^2 \psi)=1$. Then, by the definition of $w_v$, we would have $w_v(\psi)$=$(1, 0, 1)$, which is a contradiction. 
Then, $v(\varphi)$=$v(\n\psi)$=$1$, necessarily. Besides that, if 
$w_v(\psi)$ = $(0,1,0)$, then (since $\psi \in \V$), $v(\varphi)$ = $v(\n \psi)$ = $1$, from the definition of $w_v$.
\\
\noindent 2.1.1.2) Suppose $\psi$ = $\n\alpha$, with $\alpha \in \V$ (and therefore $\varphi$=$\n^2\alpha$):
if
$w_v(\psi)$ = $(1,1,0)$ then (taking into account the form of $\n$ developed in Remark \ref{M2-en-limpio})  $w_v(\alpha)$ = $(1,1,1)$ or 
$w_v(\alpha)$ = $(0,1,1)$. In both cases, since $\alpha \in \V$ and recalling definition of $w_v$, $v(\varphi)$ = $v(\n^2 \alpha)$ = $1$. 
In a similar way, if 
$w_v(\psi)$ = $(0,1,0)$ then  $w_v(\alpha)$ = $(1,0,1)$ and so $v(\varphi)$ = $v(\n^2 \alpha)$ = $1$.


\noindent 2.1.1.3) In all the other cases (that is, when $\psi \notin K_2^{\ast}$), note that $v(\varphi)$ = $v(\n \psi)$ = $-v(\psi)$. So, according Remark 
\ref{M2-en-limpio}, it only can happen that $w_v(\psi)$ = $(0,1,0)$ or $w_v(\psi)$ = $(1,0,1)$. Hence, $w_v(\psi)$ = $(0,1,0)$. By IH, $v(\psi)$ = $0$ and so 
$v(\varphi)$ = $1$. 


\noindent 2.1.2) $w_v(\varphi)$ =$(1,1,0)$. Here, $w_v(\psi)$=$(1,1,1)$ or $w_v(\psi)$=$(0,1,1)$ (from the truth-tables of ${\bf A_2}$). This implies, additionally, that 
$\psi\in \V$ and therefore $w_v(\psi)$ = $(v(\psi),v(\n \psi), v(\n^2 \psi))$. Then, in any case
$v(\varphi)$ = $v(\n \psi)$ = $1$.
%
%
%

\noindent 2.1.3) $w_v(\varphi)=(1,1,1)$. This cannot happen, from the truth-table of ${\bf A_2}$.

\noindent From all this, $w_v(\varphi) \in D_2$ implies $v(\varphi)$ = $1$. 

\noindent 2.2) Let us suppose that $v(\varphi)$ = $1$, now.

\noindent 2.2.1) If $\psi\in \V$, then $w_v(\psi)$ = $(v(\psi), 1, v(\n^2 \psi))$. So (recalling Definition \ref{matriz-a2-defi}), 
$w_v(\varphi)$ = 
$w_v(\n \psi)$ = $(1, v(\n^2 \psi), - v(\n^2 \psi))$. This implies that $w_v(\varphi) \in D_2$.

\noindent 2.2.2) Suppose now $\psi$ = $\n \alpha$, with $\alpha \in \V$. We have two possibilities here: if $v(\alpha)$ = $0$, the (since $v \in S_2$), 
$v(\psi)$ = $1$. Since $v(\varphi)$ = $v(\n \psi)$ =  $1$, in this case $w_v(\alpha)$ = $(0,1,1)$. So (checking the truth-tables of ${\bf A_2}$),
$w_v(\psi)$ = $(1,1,0)$ and then $w_v(\varphi)$ = $(1,0,1) \in D_2$. On the other hand; if $v(\alpha)$ = $1$, then (recalling that $v \in S_2$), we have the following possibilities for $w_v(\alpha)$ = $(v(\alpha), v(\psi), v(\varphi))$: $w_v(\alpha)$ = $(1,0,1)$ or $w_v(\alpha)$ = $(1,1,1)$ or $w_v(\alpha)$ = $(1,1,0)$. 
The last one contradicts $v(\varphi)$ = $1$. Then, applying the truth-table of  ${\bf A_2}$ for the other possibilities, $w_v(\varphi) \in D_2$.

%
%

\noindent 2.2.3) In any other case, $\psi \notin K_2^{\ast}$. Then, $v(\varphi)$ = $-v(\psi)$. So, $v(\varphi)$ implies $v(\psi)$ = $0$. BY IH, 
$w_v(\psi) \notin D_2$. Hence, $w_v(\psi)$ = $(0,1,0)$ or $w_v(\psi)$ = $(0,1,1)$. So, $w_v(\varphi)$ = $(1,0,1)$ or
$w_v(\varphi)$ = $(1,1,0)$, respectively. In both cases, $w_v(\varphi) \in D_2$.

So, from 2.1) and 2.2) $(\circledast)$ is valid for $\varphi$ = $\n \psi$.


\noindent \underline{Case 3}: $\varphi$ =$ \psi\supset\theta$. From Remark \ref{M2-en-limpio},  $w_v(\varphi) \in D_2$ if and only if $w_v(\psi) \notin D_2$ or $w_v(\theta) \in D_2$; if and only if (by IH) $v(\psi)$ = $1$ or $v(\theta)$ = $0$, if and only if (by Definition \ref{semantica-s2} ${\bf (4)}$), $v(\varphi)$ = $1$. So, $(\circledast)$ is also valid here. 
\hfill $\Box$

\begin{defi}\label{val2-to-biv2} 
\rm{Let $w:L(C) \rightarrow A_2$ be a $M_2$-valuation. 
Denoting, for every $\varphi \in L(C)$, $w(\varphi)$:=$(x_0^{\varphi},x_1^{\varphi}, x_2^{\varphi})$, with $x_i^{\varphi} \in \{0,1\}, \ 0\leq i\leq 2$, we define
\\
$v_w: L(C) \rightarrow \{0,1\}$ in a recursive way, by:
\\
\noindent $-$ If $\varphi \in \V$, then $v_w(\varphi)$ := $x_0^{\varphi}$.

\noindent $-$ If $\varphi$ = $\n \alpha$, $\alpha \in \V$, then $v_w(\varphi)$ := $x_1^{\alpha}$.

\noindent $-$ If $\varphi$ = $\n^2 \alpha$, $\alpha \in \V$, then $v_w(\varphi)$ := $x_2^{\alpha}$.

\noindent  $-$ If $\varphi$ = $\n \psi$, with $\psi \notin K_2^{\ast}$, then $v_w(\varphi)$ := $- v_w(\psi)$.

\noindent $-$ If $\varphi$ = $\psi \supset \theta$, then $v_w(\varphi)$ := $v_w(\psi) \rar v_w(\theta)$.
}
\end{defi}

\begin{prop}\label{biv-2-bien-definida}\rm{The map $v_w$ defined above is  a ${Ciu}^2$-bivaluation.}
\end{prop}

\dm Obviously $v_w$ is a well defined map from $L(C)$ to $2$. Let us see that it verifies the conditions recquired to be a 
${{Ciu}^2}$-bivaluation, cf. Definition \ref{semantica-s2}.
\\
\noindent ${\bf (1)}$ $v_w(\n\varphi)=0$ implies $v_w(\varphi)=1$. Let us suppose $v_w(\n\varphi)=0$, and let us analyze the following possibilities:
\\
\noindent \underline {Case 1:} $\varphi\in \V$. Since $v_w(\n \varphi)$ = $x_1^{\varphi}$, then $w(\varphi)$ = $(1,0,1)$. Hence, $v_w(\varphi)$ = $x_0^{\varphi}$ = $1$.
\\
\noindent \underline {Case 2:} $\varphi$ = $\n \alpha$, $\alpha \in \V$. Then, $v_w(\n \varphi)$ = $v_w(\n^2 \alpha)$ = $x_2^{\alpha}$ = $0$. Checking $A_2$, we have that $w(\alpha)$ = $(0,1,0)$ or $w(\alpha)$ = $(1,1,0)$. In both cases, $v_w(\varphi)$ = $x_1^{\alpha}$ = $1$.
\\
\noindent \underline {Case 3:} $\varphi \notin K_2^{\ast}$. Then, $0$ = $v_w(\n \varphi)$ =  $- v_w(\varphi)$, and therefore $v_w(\varphi)$ = $1$.
\\
So, ${\bf (1)}$ is valid for every $\varphi \in L(C)$.

\noindent ${\bf (2{.}2)}$  $v_w(\n^3\varphi)$ = $1$ implies $v_w(\n^2\varphi)$ = $0$. Note that $\n^2 \varphi \notin K_2^{\ast}$.
Then, $1$ = $v_w(\n^3\varphi)$ = $-v_w(\n^2\varphi)$ and so $v_w(\n^2\ \varphi)$= $0$. 
\\
\noindent ${\bf (3)}$  $v_w(\n(\varphi\supset\psi))$=$1$ implies $v_w(\varphi\rightarrow\psi)$=$0$. Again, since $\varphi \supset \psi \notin K_2^{\ast}$, proceed as in the proof of  ${\bf (2{.}2)}$ above. 
\\
\noindent ${\bf (4)}$ $v_w(\varphi\supset\psi)$ = $1$ if  and only if $v_w(\varphi)$=$0$ or $v_w(\psi)$=$1$. This is valid by similar reasons to the previous cases, considering that $v_w(\varphi \supset \psi)$ = $v_w(\varphi) \rightarrow v_w(\psi)$, since  $\varphi \supset \psi \notin K_2^{\ast}$. This concludes the proof. \hfill $\Box$

\

\begin{lema}\label{lema-clave-ciu2-2}\rm{For every $M_2$-valuation $w$, the $Ciu^2$-bivaluation $v_w$ given in Definition \ref{val2-to-biv2}  verifies, for every $\varphi \in L(C)$:
\begin{center} $(\star)$
$w(\varphi) \in D_2$ if and only if \ $v_w(\varphi)$ = $1$.
\end{center}
}
\end{lema}

\dm Let $w$ be an arbitrary $M_2$-valuation, $v_w$ defined cf. Definition \ref{val2-to-biv2}, and let $\varphi$ be in $L(C)$. We will proof 
$(\star)$ by induction on the completity of $\varphi$, again. Recall here that $D_2$ = $\{(1,0,1); (1,1,0); (1,1,1)\}$.
\\
\noindent \underline {Case 1:} $\varphi \in\mathcal{V}$. Then, $w(\varphi)$ = $(a_0^{\varphi},a_1^{\varphi},a_2^{\varphi}) \in D_2$ if and only if 
$a_0^{\varphi}$ = $1$ = $v_w(\varphi)$, cf. Def. \ref{val2-to-biv2}.
\\
\noindent \underline {Case 2:} $\varphi$ = $\n \psi$. Then, $w(\varphi)$ = $(a_0^{\varphi},a_1^{\varphi},a_2^{\varphi})$:= $(a_1^{\psi}, a_2^{\psi}, - a_2^{\psi})$, and therefore
$v_w(\varphi)$ = $a_1^{\psi}$. Now, let us consider the following sub-cases:
\\
\noindent 2.1) $\psi \in\mathcal{V}$. Then, 
$w(\varphi) \in D_2$ if and only if $v_w(\varphi)$ = $1$, checking the behavior of the truth-table of $\n$ for ${\bf A_2}$ developed in Remark \ref{M2-en-limpio}.
%
%
%
%
%
%
%
%
\noindent 2.2) $\varphi=\n\psi$ con $\psi\not\in\mathcal{V}$.
\\
\noindent 2.2.1) $\psi$ = $\n\alpha, \alpha \in \V$. That is, $\varphi$ = $\n^2 \alpha$ and so $v_w(\varphi)$ = $a_2^{\alpha}$, meanwhile 
$w(\varphi)$ = $(a_2^{\alpha}, - a_2^{\alpha}, a_2^{\alpha})$. Also here $w(\varphi) \in D_2$ if and only if $a_2^{\alpha}$ = $v_w(\varphi)$ = $1$ (checking the truth-table of ${\bf A_2}$ again).
%
\\
\noindent 2.2.2) $\psi \notin K_2^{\ast}$. So, $w(\psi)$ = $(1,0,1)$ or $w(\psi)$ = $(0,1,0)$, as we said before. From this, and taking into account Definition \ref{matriz-a2-defi} (a), $w(\varphi)$ = $(a_1^{\psi}, a_2^{\psi}, - a_2^{\psi})$ = $(-a_0^{\psi}, - a_1^{\psi}, - a_2^{\psi})$. This implies the following: 
$w(\varphi) \in D_2$ if and only if \\
$w(\varphi)$ = $(1,0,1)$, if and only if 
$w(\psi)$ = $(0,1,0)$, if and only if (by IH), $v_w(\psi)$ = $0$, if and only if $v_w(\varphi)$ = $1$, considering Definition \ref{val2-to-biv2}.

\noindent \underline{Case 3:} $\varphi$ = $\psi \supset \theta$. According Remark \ref{M2-en-limpio} $(d)$, $w(\varphi) \in D_2$ if and only if $w(\psi) \notin D_2$ or $w(\theta) \in D_2$, if and only if (by HI) $v_w(\psi)$ = $0$ or $v_w(\theta)$ = $1$. This last fact is equivalent to $v_w(\psi) \rar v_w(\theta)$ = $v_w(\varphi)$ = $1$, cf. Definition \ref{val2-to-biv2}. 

From all this, $(\star)$ is valid for every $\varphi \in L(C)$. \hfill $\Box$

\

From Lemmas \ref{lema-clave-ciu2-1} y \ref{lema-clave-ciu2-2} we can prove:

\begin{teor}\label{completitud-ciu-2}\rm
{For every $\Gamma \cup \{\varphi\} \subseteq L(C)$, $\Gamma \models_{S_2} \varphi$ if and only if $\Gamma \models_{M_2} \varphi$.}
\end{teor}

\dm Suppose first  $\Gamma \not\models_{M_2}\varphi$. Then, there is a $M_2$-valuation $w$ such that $w(\Gamma) \su D_2$,  $w(\varphi) \notin D_2$. Then, the 
$Ciu^2$-bivaluation $v_w$ given in Definition \ref {val2-to-biv2} verifies $v_w(\Gamma) \su \{1\}$, $v_w(\varphi)$ = $0$, cf. Lemma \ref{lema-clave-ciu2-2}. This means that $\Gamma \not\models_{S_2}\varphi$. Summarizing, $\Gamma \models_{S_2}\varphi$ implies $\Gamma \models_{M_2}\varphi$. For the other implication proceed in a similar way, using Lemma \ref{lema-clave-ciu2-1}. \hfill $\Box$

\begin{obss}\label{observaciones-posta}
\rm{The definitions and results given in this section can seem somewhat combinatorial and ``unnatural'' ones, until now. However, there is a hidden motivation behind them, that will be better explained later. So, we conclude this section commenting the underlying ideas developed here, which will be formalized and generalized in the sequel.

\

\noindent ${\bf (a)}$ Note first that, fixed $\alpha \in \V$, every bivaluation $v \in S_2$ evaluating the formulas $\n^j (\alpha)$, with $j \in \omega$,  depends only of the values  $\n^k \alpha$, with $0 \leq k \leq 2$, {\it with certain constraints}. Moreover, there are only five possibilities that such values can take. These cases correspond with the $3$-uples appearing in the definition of $A_2$, in Lemma \ref{lema-clave-ciu2-1}. Tha is, considering $\alpha \in \V$, the $Ciu^2$-bivaluations (applied to $\alpha)$ only can take on the following forms:

$v_1(\alpha)$ = $1$; $v_1(\n \alpha)$ = $1$; $v_1(\n^2 \alpha)$ = $1$.

$v_2(\alpha)$ = $1$; $v_2(\n \alpha)$ = $1$; $v_2(\n^2 \alpha)$ = $0$.

$v_3(\alpha)$ = $1$; $v_3(\n \alpha)$ = $0$; $v_3(\n^2 \alpha)$ = $1$.

$v_4(\alpha)$ = $0$; $v_4(\n \alpha)$ = $1$; $v_4(\n^2 \alpha)$ = $1$.

$v_5(\alpha)$ = $0$; $v_5(\n \alpha)$ = $1$; $v_5(\n^2 \alpha)$ = $0$.

\noindent In addition, for every $1 \leq r \leq 5$, the value $v_r(\n ^{j+1} \alpha)$ (with $j \geq 2$) is {\it univocally determined} by  
$v_r(\n^j \alpha)$, having in mind conditions ${\bf (1)}$ and ${\bf (2.n)}$ in Definition \ref{ciun-bivaluaciones}. More specifically, if $j \geq 2$, then 
$v_r(\n ^{j+1} \alpha)$ = $- v_r(\n^j \alpha)$).  Moreover, the same definition requires that
$v_r(\varphi \supset \psi)$ = 
$v_r(\varphi) \rar v_r(\psi)$, for any $\varphi$, $\psi \in L(C)$. From all this, if $\varphi$ = $\varphi(\alpha_1,\dots,\alpha_m)$, and $v \in S_2$, then  $v(\varphi)$ depends on the family $V(\varphi)$:= $\{v(\n^k \alpha_i): 1 \leq i \leq m; \, 0\leq k \leq 2\}$.

\noindent ${\bf (b)}$ On the other hand, the values of $V(\varphi)$ {\it are not completely independent from each other}. From conditions ${\bf (1)}$ and ${\bf (2.n)}$, for every formula $\n^k \alpha$, with $\alpha \in \V$, $0 \leq k \leq 1$): if $v(\n^k \alpha)$ = $0$, then $v(\n^{k+1}\alpha)$ = $1$. {\it However}, if $v(\n^k\alpha)$ = $1$, then {\it there are not restrictions for $v(\n^{k+1}\alpha)$}. All this justifies that, fixed $\alpha \in \V$, $v \in S_2$, the different possibilities for $v(\alpha)$, $v(\n \alpha)$, $v(\n^2 \alpha)$ only can be 
the five cases developed on $(a)$. So, if $\varphi$ = $\varphi(\alpha_1,\dots,\alpha_m)$, then $|V(\varphi)|$ = $5^m$.
All this anticipates, in an informal way, the reason because $|A_2|$ = $5$. 

\noindent ${\bf (c)}$ Besides that (and fixed $\alpha \in \V$, and $v \in S_2$), we can codify the family $\{v(\alpha)$, $v(\n \alpha), v(\n^2 \alpha)\}$ by means of a certain specific $3$-uple $(x_0, x_1, x_2)$ = $(v_r(\alpha), v_r(\n \alpha), v_r(\n^2 \alpha)) \in 2^3$ (which will be called an ``initial sequence'', from now on). So, every truth-value $\vec{x} \in A_2$ can be  naturally represented as some of the $3$-uples of $2^3$ codifying $v_1$-$v_5$, because $\vec{x}$ should be associated (in an unique way) to a specific initial sequence.

\noindent ${\bf (d)}$ Now, once the set $A_2$ of truth-values has been established, the operations $\supset$ and $\n$ in ${\bf A_2}$ and the set $D_2$ are defined obbeying the following criterium:
it should be a correspondence between the bivaluations $v \in S_2$ and the matrix valuations $w:L(C) \to A_2$, in such a way that any  $v \in S_2$ could be interpreted (in the context of the operations in $M_2$) as being the truth-value 
$w_v(\varphi)$ = $(v(\varphi), v(\n \varphi), v(\n^2 \varphi)) \in A_2$. This idea motivated the definition of $\n$, $\supset$ and also of $D_2$. Moreover,  $D_2$ indicates the following: given $w_v: L(C) \to 2$, $w_v(\varphi) \in D_2$ if and only if the bivaluation $v \in S_2$, underlying to $w$, verifies 
$v(\varphi)$ = $1$.

\

Of course, beyond the previous motivations, it should checked that, {\it actually}, the correspondence suggested behaves in the right way. Indeed, this is proved in Lemmas \ref{lema-clave-ciu2-1} and \ref{lema-clave-ciu2-2}, which implies Theorem \ref{completitud-ciu-2}. In addition, the essential point here is all this process (focused on $Ciu^2$ here) can be given in a general way for every $n \in \omega$. We will develop all this in the next section.

}
\end{obss}

\section{General Matrix Semantics for the $Ciu^n$-logics}

Until now, we have shown the matrix characterization of $Ciu^0$ ( = $CL$), 
$Ciu^1$ (= $P^1$) and $Ciu^2$. We will show here that all this can be generalized: every logic $Ciu^n$, with $n \geq 3$, can be naturally associated to a finite matrix (that we will call $M_n$). For that, we will take into account, as starting point, some comments of Remarks \ref{observaciones-posta}. Note first that the truth-values of the sets $A_n$ (the supports of  
$M_n$, to be defined) should be identified with $(n+1)$-tuples $\vec{x}$:=$(x_0,\dots, x_n)$. That tuples will codify certain {\it initial sequences}, such as it was 
done for $C^2$ (and, in a hidden way, for $C^0$ and for $C^1$, as we shall see).   
Besides that, the operations $\n$ and $\supset$ in every matrix $M_n$ should be defined as being a generalization of 
Definition \ref{matriz-a2-defi}. Bearing this in mind, we will focus our attention on the connections  between the conditions asked in  Definition \ref{ciun-bivaluaciones} and $|A_n|$, the cardinality of the set $A_n$. For that, we remark again that such cardinality is related with {Fibonacci{'}s sequence}. To prove this fact we will use the following definitions and results:

 \begin{defi}\label{fibo-defi}\rm{Consider $\mathbb{N}$:=$\omega\setminus\{0\}$ = $\{1,2,3,\dots\}$. The {\it (standard) Fibonacci{'}s sequence} is given by the map $Fb: \mathbb{N} \to \mathbb{N}$, recursively defined on  $\mathbb{N}$ as usual:
$Fb(1)$ = $1$; $Fb(2)$ = $1$; for $k \geq 3$, $Fb(k)$ = $Fb(k-1) + Fb(k-2)$. The first terms of this sequence are: 
\bc $1$, $1$, $2$, $3$, $5$, $8$, $13...$ \ec
}
\end{defi}

It is well-know that Fibonacci{'}s sequence can be obtained in several ways. In particular, it can be defined by means of its  {\it Binary Expansion} (see \cite{lot:02}), which can be formalized as  follows:

\begin{defi}\label{binaria-fibonacci}\rm{Let $L_{Bin}$ the {\it binary language}, consisting of all the finite words determined by the set $2$:= $\{0,1\}$ (interpreted in this context as the {\it alphabet} that generates $L_{Bin}$), and let $\sigma: 2 \to L_{Bin}$ be the substitution defined as follows: 
 $\sigma(0)$ = $1$; $\sigma(1)$ = $1 0$. The map $\sigma$ can be extended to every word 
$a_1 a_2 \dots a_n \in L_{Bin}$ as being the following concatenation: $\sigma (a_1, \dots, a_n)$:= $\sigma(a_1)\dots\sigma(a_n)$. Here, the expression ``$\sigma(a_1)\dots\sigma(a_n)$'' obviously means ``$\sigma(a_1)$ concatenated with $\sigma(a_2)$ concatenated with  ... with $\sigma(a_n)$''. Let us define, in addition,   
{\it Fibonacci{'}s binary expansion} as being the map  $W: \mathbb{N} \to L_{Bin}$, recursively by  $W(1)$ := $0$; $W(k+1)$ := $\sigma(W(k))$.\footnote{The substitution given in \cite{lot:02} is defined by $\sigma(0)$ = $0 1$; $\sigma(1)$ = $0$. Our adaptation determines words with the same lenght of the original approach, obviously.}}
\end{defi}

Note that the lenght of every word $W(n)$ is, precisely, $Fb(n)$. An usual graphic ``tree-description'' of $W(n)$ can be exemplified as follows:

\begin{exem}\rm{Starting from the word $W(1)$ = $0$, the words $W(1)$ - $W(5)$ can be associated to each {\it level} of the tree displayed on Figure 1 (with branches $b_1$- $b_5$ = $Fb(5)$):

\
\tiny
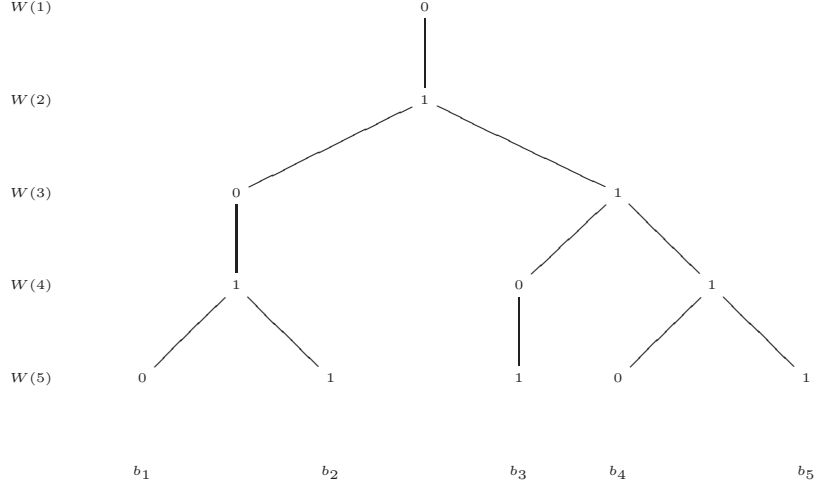
$$\xymatrix{
W(1)& & & & 0 \ar@_{-}[d] & & & & \\
W(2) & & & & 1 \ar@_{-}[drr] \ar@_{-}[lld] & & & &  \\
W(3) & & 0\ar@_{-}[d]&  & & & 1  \ar@_{-}[dr] \ar@_{-}[ld] &  \\
W(4) & & 1 \ar@_{-}[dr] \ar@_{-}[ld] &  &  & 0 \ar@_{-}[d] &  & 1 \ar@_{-}[dr] \ar@_{-}[ld] \\
W(5) & 0 &  & 1 &  & 1 & 0 & & 1 \\
& b_1 & & b_2 & & b_3 & b_4 & & b_5\\
}
$$
}
\normalsize
\end{exem}

Following the approach of the previous example,  we can consider every word $W(n)$ as a level of a tree that describes, for every $k \leq n$, ``the growth'' of every word $W(k)$ with respect to $W(k-1)$. Note now that  {\it every branch } $b_{j}$ (with $1 \leq j \leq Fb(n)$) of this tree can be interpreted (at least) by two different ways:

\

\noindent $\bullet$ If it is fixed $\alpha \in \V$, every branch $b_1,\dots, b_{Fb(n)}$ codifies all the different possible assignations of $0$ and $1$ to the formulas $\alpha$, $\n \alpha$, $\n^{Fb(n)-3} \alpha$ (the already mentioned ``initial sequences'' in Remarks \ref{observaciones-posta}). For instance, in the tree developed on Figure 1, every branch can be codified (omitting the irrelevant information of the levels $W(1)$  and $W(2)$) as follows:

$b_1$ = $(1,1,1)$

$b_2$ = $(1,1,0)$

$b_3$ = $(1,0,1)$

$b_4$ = $(0,1,1)$

$b_5$ = $(0,1,0)$


\noindent So, every branch can be interpreted as a possible initial sequence (for the logic $Ciu^2$, in this example). 

\

\noindent $\bullet$ In addition, the $3$-uples that identify the branches $b_1 - b_5$ are intended as being the truth-values of $A_2$, and so $|A_2|$ = $5$ = $Fb(5)$. 

\noindent Note now that, in the case of $Ciu^0$ (= $CL$), $|A_0|$ = $|\{0,1\}|$ = $2$ = $Fb(3)$. 
Moreover, this relation is also valid for $Ciu^1$ (i.e. Sette{'}s logic $P^1$), whose support has  three truth-values. In fact, in \cite{ram:fer:09} it is shown that the set of the truth-values of $P^1$ can be identified with the set 
$A_1$ = $\{(1,1); (1,0); (0,1) \}$ \footnote{By the way, the general structure of the truth-values of the sets $A_n$ that we will propose in this section is based on the identification given in the mentioned paper, which is developed in a deeper way in \cite{fer:mur:18}.}. In other words, $|A_1|$ = $Fb(4)$.

\

The relation between the initial sequences referred to $Ciu^n$ (which will be the truth-values of $A_n$) and Fibonacci{'}s Binary Expansion is very natural, indeed. It is based on the fact that, if $\varphi$ = $\n^k \alpha$, $0 \leq k < n$ and $v \in S_n$, it is recquired that  $v(\varphi)$ = $0$ implies $v(\n \varphi)$ = $1$, but {\it there are no restrictions} when $v(\varphi)$ = $1$ (this determines an obvious bifurcation). 
And this behavior can be codified by means of 
$\sigma: 2 \to L_{Bin}$, in the context of the word $W(n)$!

\

With this idea in mind, it is very natural to conjecture that, if  $M_n$ = $({\bf A_n},D_n)$ is the matrix semantics associated to $Ciu^n$, then 
$|A_n|$ = $Fb(n+3)$. Unifying this with the comments of Remarks \ref {observaciones-posta}, it is possible to define $M_n$ in a natural way. We will develop the technical aspects of all these comments  in the sequel. We remark here that many definitions and proofs of the results to be shown from now on {\it are different} from the similar ones given on Section \ref{semantica-para-ciu2}, because the general analysis for every logic $Ciu^n$ requires a more abstract approach.

\begin{defi}\label{soportes-ciu-n-general}\rm{Let us consider the set  $2$ = $\{0,1\}$: 

\noindent $(a)$ The set $A_n \su 2^{n+1}$ is defined (by recursion on $n \in \omega$) by:

 \noindent $-$ $A_0$ := $\{0,1\}$ (= $2$).

\noindent $-$ From $A_{n-1}$, it is defined  $A_{n}$ as follows: 

$A_{n}$ := $\{\vec{x} \in 2^{n+1}: \vec{x} \textrm{ verifies condition $\diamond$}\}$. 

\noindent Here, condition $\diamond$ is:
$\vec{x}$ = $(x_0,\dots, x_{n-1},x_n)$ verifies $\diamond$ if and only if:

\noindent $(\diamond{.}1)$ $(x_0,\dots, x_{n-1}) \in A_{n-1}$

\noindent $(\diamond{.}2)$ If $x_{n-1}$ = $0$, then $x_{n}$ = $1$

\noindent  (There are not restrictions for $x_{n}$, if $x_{n-1}$ = $1$).

\noindent $(b)$ In addition, for every $n \geq 0$, we define $D_n$:= $\{ \vec{x}\in A_n: x_0 = 1 \}$.
}
\end{defi}

\begin{obs}\label{An-Fibo}\rm{
From this recursive definition it can be easily seen that, for every $n \geq 0$, $|A_n|$ = $Fb(n+3)$. Moreover, it holds that $|D_{n+1}|$ = $|A_n|$. On the other hand, note that Definition \ref{soportes-ciu-n-general} is equivalent to the following  non-recursive expression, which is a generalization of Definition \ref{matriz-a2-defi}:
}
\end{obs}

\begin{prop} \label{def-equivalentes-An}
\rm{For every $n \in \omega$, the set $A_n$ can be characterized as follows:
$A_n$ = 
$\{\vec{x} = (x_0,\dots,x_{n-1},x_n)\in 2^{n+1}: 
\textrm{$x_k = 0$ implies $x_{k+1} = 1$, 
$0 \leq k 
\leq n-1$}\}$.
} 
\end{prop}

\begin{defi}\label{matrices-ciu-n-general}\rm{For every $n \geq 0$, we define the $C$-matrix $M_n$:= $({\bf A_n},D_{n})$, where the support of 
${\bf A_n}$ is $A_n$, and the operations in $M_n$ are defined as follows (for $\vec{x}$ = $(x_0,\dots,x_n)$, 
$\vec{y}$ = $(y_0,\dots, y_n) \in A_n$):

\noindent $(a)$ $\n \vec{x}$:= $(x_1,x_2,\dots,x_n, - x_n)$; in the case of $A_0$ this operation must be understood in this way: $\n \vec{x}$ = $\n(a_0)$:=$-a_0$.

\noindent $(b)$ $\vec{x} \supset \vec{y}$:=$\vec{z}$ = $(z_0,\dots,z_n)$, where: 
$z_0$ = $x_0 \rar y_0$; $z_{k}$:= $-z_{k-1}$, for every $1 \leq k \leq n$.

}
\end{defi}

\begin{exems}
\rm{

\

\

\noindent $\bullet$ When $n$ = $0$, $A_0$ = $\{0,1\}$; $D_0$ = $\{1\}$, $\n \vec{x}$  = $\n x_0$:=$- x_0$ and $\vec{x} \supset \vec{y}$ = $x_0 \supset y_0$:= $x_0 \rar y_0$.

\noindent $\bullet$ If $n$ = $1$, then $A_1$ = $\{(0,1); (1,0); (1,1) \}$, $D_1$ = $\{(1,0); (1,1)\}$. This case coincides with the matrix for $P^1$ as presented in \cite{ram:fer:09}, as it is was commented.

\noindent $\bullet$ The matrix $M_2$ (developed according Definition \ref{matrices-ciu-n-general}) is, actually, the same matrix given in Remark \ref{M2-en-limpio}.

\noindent $\bullet$ If $n$ = $3$, $A_3$ = $\{(0,1,0,1); (0,1,1,0);, (0,1,1,1); (1,0,1,0); (1,0,1,1); (1,1,0,1);$ 

\noindent $(1,1,1,0); (1,1,1,1) \}$; $D_3$ = $\{(1,0,1,0); (1,0,1,1); (1,1,0,1); (1,1,1,0); (1,1,1,1)\}$. We omit here the behavior of $\supset$ and of $\n$.
}
\end{exems}

From now on, we will often use the set $K_n^{\ast}$ (remember it from Definition \ref{k-n}). As a starting point, we will relate it with the $M_n$-valuations and with the truth values of $A_n$ defined above.

\begin{prop}\label{propiedades-grales-kn-ast}\rm{For every $n \in \omega$, the set $K_n^{\ast}$ verifies:
%
%
%

\noindent $(a)$ For every $\varphi \notin K_n^{\ast}$, for every $M_n$-valuation $w$, it only can happen one of these possibilities:

\noindent $(a{.}1)$ $w(\varphi)$ = $(x_0^{\varphi},\dots,x_n^{\varphi})$ where, for every $0 \leq i \leq n$: 
$x_i^{\varphi}$= $\left\{
\begin{array}{ll}
0 &  \textrm{if $i$ is even}\\
1 & \textrm{if $i$ is odd}
\end{array}
\right.$

\

\noindent $(a{.}2)$ $w(\varphi)$ = $(x_0^{\varphi},\dots,x_n^{\varphi})$ where, for every $0 \leq i \leq n$: 
$x_i^{\varphi}$= $\left\{
\begin{array}{ll}
0 &  \textrm{if $i$ is odd}\\
1 & \textrm{if $i$ is even}
\end{array}\right.$

\

Moreover, from $(a)$ it holds:

\noindent $(b)$ For every pair $\varphi$, $\psi \notin K_n^{\ast}$, if $w$ is an $M_n$-valuation where $w(\varphi)$ = $(x_0^{\varphi},\dots,x_n^{\varphi})$,
 $w(\psi)$ = $(x_0^{\psi},\dots,x_n^{\psi})$ :

$(b{.}1)$ For every $0 \leq i \leq n-1$, $x_i^{\varphi}$ = $- x_{i+1}^{\varphi}$.

$(b{.}2)$ $w(\n \varphi)$ = $(- x_0^{\varphi}, \dots, - x_n^{\varphi})$.

$(b{.}3)$ $w(\varphi \supset \psi)$ = $(x_0^{\varphi} \rar x_0^{\psi},x_1^{\varphi} \rar x_1^{\psi},\dots,x_n^{\varphi}\rar x_n^{\psi})$.
}
\end{prop}

From now on we will show some results connecting bivaluations of $S_n$ with $M_n$-valuations.

\begin{defi}\label{vw-defi-general}\rm{ Let $v:L(C)\to 2$ be in $S_n$. We define the $M_n$-valuation $w_v:L(C) \to A_n$ in this way, for every $\alpha \in \V$: \\
$w_v(\alpha)$ := $(a_{0}^{\alpha}, a_{1}^{\alpha}, ... , a_{n}^{\alpha})$:= $(v(\alpha), v(\n \alpha), ... ,v(\n^{n}\alpha))$. We extend 
\\
$w_v: L(C) \to A_n$ homomorphically for every $\varphi \in L(C)$, following Definition \ref{matrices-ciu-n-general}. That is:
\\
\noindent $(i)$ \ If $w_{v}(\psi)$ = $(a_{0}^{\psi}, a_{1}^{\psi}, ... , a_{n}^{\psi})$ then $w_{v}(\n\psi)$ = $\n(w_v(\psi))$:=$(a_{1}^{\psi}, ... , a_{n}^{\psi}, -a_{n}^{\psi})$.
\\
\noindent $(ii)$ \ If $w_{v}(\psi)=(a_{0}^{\psi}, a_{1}^{\psi}, ... , a_{n}^{\psi})$ and $w_{v}(\theta)=(a_{0}^{\theta}, a_{1}^{\theta}, ... , a_{n}^{\theta})$ then $w_{v}(\psi\supset\theta)$=

\noindent = $w_v(\psi)\supset w_v(\theta)$:=
$(a_{0}^{\psi}\to a_{0}^{\theta}, -(a_{0}^{\psi}\to a_{0}^{\theta}), a_{0}^{\psi}\to a_{0}^{\theta}, -(a_{0}^{\psi}\to a_{0}^{\theta}),...)$.
}
\end{defi}

\begin{prop}\label{wv-bien-definida}\rm{ For every $n \in \omega$, $\varphi\in L(C)$, it holds that $w_{v}(\varphi)\in A_n$.}
\end{prop}

\dm We will use the characterization of the sets $A_n$ given in Definition \ref{def-equivalentes-An}. Suppose $n \in \omega$; if $n$ = $0$ then $A_0$ = $2$ and our claims holds trivially. For $n \geq 1$ we will use induction on the complexity of the formulas of $L(C)$:

\noindent \underline{Case 1}: $\varphi \in \V$. Then, $w_v(\varphi)$ = $(v(\varphi),\dots,v(\n^n \varphi)) \in A_{n}$, since $v \in S_n$.

\noindent \underline{Case 2}: $\varphi$ =$\n \psi$. Here, $w_v(\varphi)$ = $(a_0^{\varphi},a_1^{\varphi},\dots,a_n^{\varphi})$ = $(a_1^{\psi},\dots,a_n^{\psi},- a_n^{\psi})$. Suppose $a_i^{\varphi}$ = $0$. First, if $i<n-1$ then $a_i^{\varphi}$ = $a_{i+1}^{\psi}$. So, by (IH), $1$ = $a_{i+2}^{\psi}$ = $a_{i+1}^{\varphi}$. On the other hand, if $i$ = $n-1$ then $a_{i+1}^{\varphi}$ = $a_{n}^{\varphi}$ = $- a_{n-1}^{\varphi}$ = $1$. 

\noindent \underline{Case 3}: $\varphi$ = $\psi \rar \theta$. Let us consider $w_v(\varphi)$ = $(a_0^{\varphi},\dots,a_n^{\varphi})$: by the definition of $w_v$ we have that, for every  $0 \leq i \leq n-1$, $a_{i+1}^{\varphi}$ = $- a_i^{\varphi}$. Hence, if $a_i^{\varphi}$ = $0$, then  $a_{i+1}^{\varphi}$ = $1$ and so $w_v(\varphi) \in A_n$ here, too. \hfill $\Box$

\begin{prop}\label{val-copia-bival}\rm{ Let $n$ be in $\omega$, $n$ fixed. Consider $v \in S_n$ and the $M_n$- valuation $w_v:L(C)\to A_n$.  For every $\varphi \in L(C)$ 
it is verified 
$$(\star): w_{v}(\varphi)=(v(\varphi), v(\n\varphi), ..., v(\n^{n}\varphi)).$$
}
\end{prop}

\dm Let us fix $n \in \omega$, $v \in S_n$, the $M_n$-valuation $w_v$ and let us prove $(\star)$ by structural induction:
\\
\noindent \underline{Case 1:} When $\varphi\in\mathcal{V}$, $(\star)$ is valid because Definition \ref{vw-defi-general}.
\\
\noindent \underline{Case 2:} If $\varphi$ = $\n\psi$, then:
\\
\noindent 2.1) $\psi\in K_n^{\ast}$. That is, $\psi=\n^p\alpha$ with 
$0\leq p< n$ and $\alpha\in\mathcal{V}$. Then: \\
$w_v(\varphi)$ = $w_v(\n\psi)$ = $\n w_v(\psi)$ = $\n (v(\psi),v(\n\psi),...,v(\n^n\psi))$ =\\
$(v(\n\psi),...,v(\n^n\psi),-v(\n^n\psi))$ = $(v(\varphi),...,v(\n^{n-1}\varphi),-v(\n^{n+p}\alpha))$=
\\
$(v(\varphi),...,v(\n^{n-1}\varphi),v(\n^{n+1+p}\alpha))$ = $(v(\varphi),...,v(\n^{n-1}\varphi),v(\n^n(\n(\n^p\alpha))))$ =
\\
$(v(\varphi),...,v(\n^{n-1}\varphi),v(\n^n\varphi))$.
\\
\noindent 2.2)  $\psi\notin K_n^{\ast}$ (and then $\varphi \notin K_n^{\ast}$). Therefore, using (IH), $w_v(\varphi)$ = $w_v(\n\psi)$ = $\n w_v(\psi)$ = $\n (v(\psi),v(\n\psi),...,v(\n^n\psi))$ = 
$(v(\n\psi),...,v(\n^n\psi),-v(\n^n\psi))$ = 
\\
$(v(\varphi),...,v(\n^{n-1}\varphi),v(\n^{n+1}\psi))$ = $(v(\varphi),...,v(\n^{n-1}\varphi),v(\n^n\varphi))$.
\\
\noindent \underline {Case 3:} When $\varphi$ = $\psi\supset\theta$. By (IH), $w_{v}(\psi)$ = $(v(\psi), v(\n\psi), ...,v(\n^{n}\psi))$ and 
$w_{v}(\theta)$ = $(v(\theta), v(\n\theta), ...,v(\n^{n}\theta))$. Therefore, we have $w_{v}(\varphi)$ = $w_{v}(\psi\supset\theta)$ = $w_{v}(\psi)\supset w_{v}(\theta)$
= $(v(\psi), v(\n\psi), ...,v(\n^{n}\psi))\supset(v(\theta), v(\n\theta), ...,v(\n^{n}\theta))$ = 

\noindent $(v(\psi)\rar v(\theta),-(v(\psi)\rar v(\theta)), -(-(v(\psi)\rar v(\theta))), ...)$=$(c_0, c_1, c_2, ...)$.
\\
Now, since $v \in S_n$ and $\varphi$ = $\psi \rar \theta \notin K_n^{\ast}$, it holds that $c_0$ = $v(\varphi)$. Moreover, note that for every  $0\leq k \leq n$, 
$v(\n^{k}\varphi)$ = $\underbrace{- \dots(-}_{k \, \textrm{times}} v(\varphi))$ = $c_k$ (from Proposition \ref{propiedades-grales-kn-ast} $(b)$). So,
$w_{v}(\varphi)$ = $(v(\varphi), v(\n\varphi), ..., v(\n^{n}\varphi))$. That is, $(\star)$ is also valid in this case. \hfill $\Box$

\begin{lema}\label{gral-fundam-1}\rm{For every
$v \in S_n$, the $M_n$-valuation $w_{v}:L(C)\to A_n$ verifies the following: for every 
$\varphi\in L(C)$, $w_{v}(\varphi)\in D_n$ if and only if $v(\varphi)$=$1$}. 
\end{lema}

\dm Let $\varphi$ be in $L(C)$ and $w_v(\varphi)$ = $(a_0^{\varphi},\dots,a_n^{\varphi})$. By Definition \ref{soportes-ciu-n-general} $(b)$, $w_w(\alpha)\in D_n$ if and only if $a_0^{\varphi}$ = $1$, if and only if  $v(\varphi)$ = $1$ (by Proposition \ref{val-copia-bival}). \hfill $\Box$

\noindent \begin{defi}
\label{nbival}\rm{Let $w:L(C)\to A_n$ be an $M_n$-valuation. We define \\
    $v_{w}:L(C)\to 2$ as follows (with the convention $w(\varphi)=(a_{0}^{\varphi}, a_{1}^{\varphi}, ...,a_{n}^{\varphi})$): 
\\
\noindent $\bullet$ If $\varphi\in\mathcal{V}$, then $v_{w}(\varphi)$:=$a_{0}^{\varphi}$. 
\\
\noindent  $\bullet$ If $\varphi$ = $\n \psi$, then:

If $\psi \in K_n^{\ast}$, 
then 
$v_w(\varphi)$ := $a_{k+1}^{\alpha}$ 

(by the way; note here that $\varphi$ = $\n^{k+1}\alpha$, $1 \leq k+1 \leq n$).

If $\psi \notin K_n^{\ast}$, then $v_w(\varphi)$:= $- v_w(\psi)$.
\\
\noindent $\bullet$ If $\varphi$= $\psi\supset\theta$ then $v_w(\varphi)=v_w(\psi)\to v_w(\theta)$.
}
\end{defi}

\begin{prop}\label{nbival-biendef}
\rm{$v_w \in S_n$.}
\end{prop}

\dm $v_w:L(C) \to 2$is a well-defined function, obviously. Now, given $\varphi , \psi\in L(C)$, it holds:
\\
\noindent {\bf (1)} 
{If  $v_w(\n\varphi)=0 $ then  $v_w(\varphi)=1$.} In fact:
Let us suppose first that $\varphi \in K_n^{\ast}$ (i.e. $\varphi$ = $\n^k \alpha$, $\alpha \in \V$, $0 \leq k < n$). Considering 
$w(\alpha)$ = $(a_{0}^{\alpha},\dots,a_n^{\alpha})$ and Definition \ref{nbival}, we have that
$v_w(\n \varphi)$ = $a_{k+1}^{\alpha}$. Now: from Definition \ref{soportes-ciu-n-general} $(\diamond{.}2)$,
it must hold that $a_k^{\alpha}$ = $1$. In addition, note that $v_w(\varphi)$ = $a_k^{\alpha}$. Therefore, $v_w(\varphi)$ = $1$. On the other hand, if 
$\varphi \notin K_n^{\ast}$, then $v_w(\n \varphi)$ = $- v_w(\varphi)$, and so our claim is also valid here.
\\
\noindent {\bf (2)} 
{If $v_w(\n^{n+1}\varphi)=1$ then $v_w(\n^{n}\varphi)=0$.} Note here that $\n^{n}\varphi \notin K_n^{\ast}$. So,
$v_w(\n^{n+1}\varphi)$ = $-v_w(\n^n\varphi)$, cf. Definition \ref{nbival}), and therefore this item is trivially valid.
\\
\noindent {\bf (3)} 
{If $v_w(\n(\varphi\rar\psi))=1$ then $v_w(\varphi\rar\psi))=0$.} Here $1$ = $v_w(\n(\varphi\rar\psi))$= $-v_w(\varphi\rar\psi)$ and so $v_w(\varphi\rar\psi)$ = $0$.
\\
\noindent {\bf (4)} 
{$v_w(\varphi\rar\psi)=1$ if and only if $v_w(\varphi)=0$ or $v_w(\psi)=1$.} Indeed,
$v_w(\varphi\rar\psi)=1$ if and only if $v_w(\varphi)\rar v_w(\psi)=1$ if and only if $v_w(\varphi)=0$ or $v_w(\psi)=1$. This concludes the proof.\hfill $\Box$

\begin{prop}\label{acero}\rm{For every $\varphi \in L(C)$, for every $n \in \omega$, for every $M_n$-valuation $w: L(C) \to A_n$ with 
$w(\varphi)$ = $(a_0^{\varphi},\dots,a_n^{\varphi})$, it holds that $v_w(\varphi)$ = $a_0^{\varphi}$. }
\end{prop}

\dm By induction on the completity of $\varphi$:

\noindent \underline{Case 1:} $\varphi\in \V$. Our claim holds from the definition of $v_w$.
 
\noindent \underline{Case 2:} $\varphi$ = $\n \psi$.

\noindent $2{.}1)$ $\psi \in K_n^{\ast}$. 
So, 
$v_w(\varphi)$ = $a_{k+1}^{\alpha}$, by the definition of $v_w$). Besides that, from Definition \ref{matrices-ciu-n-general}, 
$w(\varphi)$ = $\n^{k+1} w(\alpha)$ = $(a_{k+1}^{\alpha},\dots,a_n^{\alpha}, - a_n^{\alpha},\dots)$. Then, $v_w(\varphi)$ = $a_0^{\varphi}$.

\noindent $2{.}2)$ $\psi \notin K_n^{\ast}$. By Proposition \ref{propiedades-grales-kn-ast}, 
for $w(\psi)$ = $(a_0^{\psi},a_1^{\psi},\dots,a_n^{\psi})$ it holds that $a_0^{\psi}$ = $-a_1^{\psi}$.
In addition, $v_w(\varphi)$ = $- v_w(\psi)$ = $- a_0^{\psi}$, by (IH). So, $v_w(\varphi)$ = $- (- a_1^{\psi})$ = $a_1^{\psi}$ = $a_0^{\varphi}$.   
  
\noindent \underline{Case 3:} $\varphi$ = $\psi \supset \theta$.  Here, considering $w(\varphi)$ and (IH),
\\
$v_w(\varphi)$ = $v_w(\psi \supset \theta)$ = $v_w(\psi) \rar v_w(\theta)$ = $a_0^{\psi} \rar a_0^{\theta}$ = $a_0^{\varphi}$. 
This concludes the proof. \hfill $\Box$

\begin{lema}\label{gral-fundam-2}\rm{For every $M_n$-bivaluation $w$, for every $\varphi \in L(C)$, $v_w(\varphi)$ = $1$ if and only if $w(\varphi) \in D_n$.}
\end{lema}

\dm Considering $w(\varphi)$ = $(a_0^{\varphi},\dots,a_n^{\varphi})$ and the result above, we have that $v_w(\varphi)$ = $1$ if and only if $a_0^{\varphi}$ = $1$, if and only if $w(\varphi)\in D_n$. \hfill $\Box$

\begin{teor}\label{EQ}\rm{Given $n \in \omega$, $\Gamma \cup \{\varphi\}\su L(C)$, it holds that $\Gamma \models_{S_n} \varphi$ if and only if $\Gamma \models_{M_n} \varphi$.}
\end{teor}

\dm Adapt Theorem \ref{completitud-ciu-2} for any $n \in \omega$, applying Lemma \ref{gral-fundam-1} and Lemma \ref{gral-fundam-2}. \hfill $\Box$

\section {Conclusions}

We have proved in this paper that each logic of the family $\textsf{Ciu}^n$:=$\{Ciu^n\}_{n \in \omega}$ admits finite matrix representation. Moreover, the cardinality of (the support of) such matrices grows up according Fibonacci{'}s sequence $Fib(n)$. This last fact is not a trivial one, since in many logics with finite matrix semantics (as the already mentioned hierarchy $\{I^0 P^n\}$ of \cite{fer:con:03}), that growth is lineal. More specifically in this case: 
$|I^0 P^n|$ = $n+2$. Now, in the $Ciu^n$-logics, the cardinality is directly related with the behavior of the bivaluations that defined them. This suggest an interesting topic of future work: the study of the relation between the definition of  an arbitrary bivaluated logic $\lc$, with its finite matrix characterization, if it exists. An example of this kind of situation is the following: as it was mentioned, we have that da Costa{'}s logics $C_n$ have a bivaluated semantics (underlying in the quasi-matrices provided in \cite{dac:alv:77}). Despite this fact, these logics {\it are not representable by finite matrices}. In which way the form of the bivaluations of the $C_n$-logics is related with the latter result?

Continuing with this approach, it is our aim the study of an alternative paraconsistent hierarchy (also defined by J. Ciuciura), which we will call here as
$\textsf{Ciu}^{\ast n}$:=$\{Ciu^{\ast n}\}_{n \in \omega}$. Sintetically, the main difference between the  $Ciu^n$-logics and the $Ciu^{\ast n}$-logics is that in the latter {\it it is not required} that every bivaluation $v$ behaves homomorphically w.r.t. the implication connective $\supset$. That is, conditions ${\bf (3)}$ and ${\bf (4)}$ of Definition \ref{ciun-bivaluaciones} are not valid anymore (they are replaced by weaker versions of them, cf. \cite{ciu:20}, Definition 4.5). Actually, we conjecture that the $Ciu^{\ast}_n$-logics {\it do not admit finite matrix representation}. 
The analyisis of these topics will be developed in future studies relating bivaluated semantics with matrix ones, in a general way.

Finally, we wish to remark that the finite matrix representation of the $Ciu^n$-logics implies certain technical advantages.
For instance, it is possible to furnish a proof of completeness based on certain methods that usually are applied to $n$-valued logics, such as the constructive technique of L. Kalm\'ar (see \cite{men:97}). In this sense, such proof (which we wish to adapt for the $Ciu^n$-logics) differs from the already given in \cite{ciu:20}, based on an adaptation of Lindenbaum{'}s Lemma. Moreover, it would be interesting to obtain a Kalm\'ar-style completeness proof for an arbitrary bivaluated logic
$\lc$, without necessity of a matrix characterization of it.

\end{document}